\documentclass[12pt]{amsart}
\usepackage{amssymb}
\newtheorem{theorem}{Theorem}[section]
\newtheorem{lemma}[theorem]{Lemma}
\newtheorem{proposition}[theorem]{Proposition}
\newtheorem{corollary}[theorem]{Corollary}
\theoremstyle{definition}

\theoremstyle{remark}
\newtheorem{remarks}[theorem]{Remarks}

\numberwithin{equation}{subsection}
\newfont{\ajb}{eufm10 at12pt}
\newfont{\aj}{eufm10 at10pt}
\newfont{\ajk}{eufm10 at8pt}
\newfont{\kh}{msbm10 at10pt}
\newfont{\khk}{msbm10 at 8pt}
\newcommand{\bea}{\begin{eqnarray*}}
\newcommand{\eea}{\end{eqnarray*}}

\begin{document}
\title[Arens regularity of bilinear mappings $\cdots$]{Arens regularity of  module actions and the second adjoint  of a derivation}{}
\author{S. Mohammadzadeh}
\address{Department of Mathematics, Ferdowsi University of Mashhad, P. O. Box 91775-1159, Mashhad, Iran} \email{somohammadzad@yahoo.com}
\author{H. R. E. Vishki}
\address{Department of Mathematics, Ferdowsi University of Mashhad, P. O. Box 91775-1159, Mashhad, Iran}
\email{vishki@ferdowsi.um.ac.ir}

\begin{abstract}
In this paper, first we give a simple criterion for the Arens regularity of a bilinear mapping on normed spaces, which applies in particular to Banach module actions and then we investigate those conditions under which the second adjoint of a   derivation into a dual Banach module is again a derivation. As a consequence  of the main result,  a  simple and direct proof for several older results is also included.

\end{abstract}
\maketitle
\section{Introduction}
In his pioneering paper, \cite {AR},  Arens  has shown that a bounded bilinear map $f: X\times Y\longrightarrow Z$ on normed spaces,  has two natural but, in general, different extensions to the bilinear maps $f^{***}$ and $f^{r***r}.$ When these extensions are equal, $f$ is said to be (Arens) regular. If the multiplication of a Banach algebra $A$ enjoys this property, then $A$ itself is called (Arens) regular.

In this paper we first provide a criterion for the regularity of a bounded bilinear map (Theorem~\ref{T1} below),  by showing that $f$ is regular if and only if $f^{****}(Z^*,  X^{**})\subseteq Y^*$; which in turn covers some older results of \cite {DPV}, \cite {AR2} on this topic (see Corollaries~\ref{C1} and ~\ref{C2} below). Then we apply the above mentioned criterion for the module actions of a Banach $A-$module $X$, which in turn  give rise  to the module actions of $A^{**}$ (equipped with each of the Arens products) on   $X^{**}$ and   $X^{***}$, in a natural way. In this direction we present Propositions ~\ref{p1}, ~\ref{p2} and ~\ref{p3}, which  generalize some results of \cite {AR1}, \cite {BP}, \cite {DPV}, and \cite U.

 For a Banach $A-$module $X$, the  second adjoint $D^{**}: A^{**}\rightarrow X^{***}$ of a  derivation $D:A\rightarrow X^*$ is trivially a linear extension of $D$. A problem which is of interest is under what  conditions $D^{**}$ is again a derivation.  Dales, Rodriguez-Palacios and Velasco in \cite {DPV} studied this problem for the special case $X=A$, and they showed that $D^{**}$ is a derivation if and only if $D^{**}(A^{**})\cdot A^{**}\subseteq A^*$,  (see \cite[Theorem 7.1]{DPV}). We extend their result for a general derivation $D:A\rightarrow X^*$ with a direct proof (see Theorem~\ref{T2} below). This theorem in turn extends some other results of \cite {DPV} and \cite{BP} for a general derivation  $D: A\rightarrow X^*$.

  For terminology and background materials we follow \cite D, as far as possible.

\section{Arens regularity of bilinear maps}
For a normed space $X$, we denote by $X^*$ the topological dual of $X$. We write $X^{**}$ for $(X^*)^*$, and so on. Throughout the paper, we usually identify a normed space with its canonical image in its second dual.

Let $X$, $Y$ and $Z$ be normed spaces and let $f :X \times Y\longrightarrow Z $  be a  bounded bilinear map. The
adjoint   $f^{*} :Z^{*} \times X
\longrightarrow Y^{*}$ of $f$ is defined by
\[\langle f^{*}(z^{*},x), y\rangle =\langle z^{*},f(x,y)\rangle \ \ \ \ ( x \in X, y \in Y , z^{*} \in Z^{*} ),\]
which is also a bounded bilinear map. By setting $f^{**} = (
f^{*})^{*} $ and continuing this way, the  mappings  $ f^{**}
:Y^{**} \times Z^{*} \longrightarrow X^{*}$, $f^{***} :X^{**}
\times Y^{**} \longrightarrow Z^{**} $ and $f^{****} :Z^{***}
\times X^{**} \longrightarrow Y^{***} $ may be defined similarly.

The mapping $f^{***}$ is the unique extension of $ f $ such that $f^{***}(\cdot
, y^{**})$ is $ w^{*}-$continuous for every $y^{**} \in Y^{**}.$
Also $f^{***}(x, \cdot )$ is $ w^{*}-$continuous for every $x \in X$.

We also denote by $f^{r}$  the flip map of $f$, that is the bounded bilinear map $f^{r}: Y
\times X \longrightarrow Z$  defined by $f^{r}(y,x)= f(x,y)\  (
x \in X, y \in Y).$ It may be raised as above  a bounded bilinear map $f^{r***r} :X^{**}
\times Y^{**} \longrightarrow Z^{**} $ which in turn  is a unique extension of $f$ such that $f^{r***r}(x^{**}, \cdot)$ is $ w^{*}-$continuous for every $x^{**} \in X^{**}$, and also $f^{r***r}(\cdot , y )$ is $ w^{*}-$continuous for every $y \in Y$. One may also easily verify that \[f^{***}(x^{**}, y^{**})=w^*-\lim_\alpha\lim_\beta f(x_\alpha, y_\beta)\] and \[f^{r***r}(x^{**}, y^{**})=w^*-\lim_\beta\lim_\alpha f(x_\alpha, y_\beta),\] where $\{x_\alpha\}$ and $\{y_\beta\}$ are nets in $X$ and $Y$ which converge to $x^{**}$ and $y^{**}$  in the $w^*-$topologies, respectively.

The map $f $ is called Arens
regular when the equality $f^{***}=f^{r***r}$ holds.

It is easy to verify that,  for the multiplication map $\pi: A\times A \longrightarrow A$ of  a Banach algebra $A$,  $\pi^{***}$ and $\pi^{r***r}$ are actually the so-called
 first and second Arens products \cite {AR},  which will be denoted by
 $\Box$ and $\lozenge$,
respectively. The Banach algebra $A$ is said to be Arens regular if the multipliction  map $\pi$ is Arens regular, or equivalently $\Box=\lozenge$ on the whole of $A^{**}$. This is the case, for example, for all $C^*-$algebras, \cite {CY}, and also for $L^1(G)$ if (and only if) $G$ is finite, \cite Y. The interested reader may refer to \cite {DL} for ample information about the Arens regularity problem on a wide variety  of Banach algebras.\\

We commence with the main theorem of this section which provide a criterion concerning to the regularity of a bounded bilinear map.

\begin{theorem}\label{T1}
For a bounded bilinear map $f :X \times Y\longrightarrow Z $
the following statements are equivalent:

$(i)$  $f$ is regular.

$(ii)$  $ f^{****}=f^{r*****r}$.

$(iii)$ $ f^{****} (Z^{*}, X^{**}) \subseteq Y^{*} $.

$(iv)$ The linear map $x\longmapsto f^*(z^*, x) :X\longrightarrow Y^* $ is weakly compact for every $z^{*} \in Z^{*}$.

\end{theorem}
\begin{proof} Let $x^{**}\in X^{**}, y^{**}\in Y^{**}, z^{***}\in Z^{***}$ and $z^*\in Z^*$ be arbitrary. If $(i)$ holds then
\bea \langle f^{****}(z^{***}, x^{**}), y^{**}\rangle &=&\langle z^{***}, f^{***}(x^{**}, y^{**})\rangle\\&=&\langle z^{***}, f^{r***r}(x^{**}, y^{**})\rangle\\&=&\langle z^{***}, f^{r***}(y^{**}, x^{**})\rangle\\&=&\langle x^{**}, f^{r****}(z^{***}, y^{**})\rangle \\&=&\langle f^{r*****}(x^{**}, z^{***}), y^{**}\rangle.\eea
Therefore  $f^{****}=f^{r*****r}$, as claimed.

 The implication $(ii)\Rightarrow (iii)$ follows trivially from the fact that \[f^{****}(z^{*}, x^{**})=f^{r*****}(x^{**}, z^{*})={f^{r*****}}_{|_{X^{**}\times Z^{*}}}(x^{**}, z^{*})=f^{r**}(x^{**}, z^{*})\in Y^{*}.\]

 That $(iii)$ is equivalent to $(iv)$ is obvious; indeed,  if we denote $f^*(z^*, \cdot )$ by $L$, then it is easy to
show that
$L^{**}=f^{****}(z^{*}, \cdot )$. Now the conclusion follows from the fact  that  $L$ is weakly compact if and only if
$L^{**}(X^{**})\subseteq Y^{*}.$

For $(iii)\Rightarrow (i)$, suppose that  $f^{****} (Z^{*},X^{**}) \subseteq Y^{*} $ and
  let $\{x_{\alpha}\}$, $\{y_{\beta}\}$ be two nets in $X$ and $Y$ that converge to $x^{**}$, $y^{**}$ in the  $w^*-$topologies, respectively; then \bea \langle f^{***}(x^{**}, y^{**}),  z^{*}\rangle&=&\langle f^{****}(z^{*}, x^{**})
, y^{**}\rangle\\&=&\lim_{\beta}\langle f^{****}( z^{*}, x^{**})
, y_{\beta}\rangle\\&=&\lim_{\beta} \langle f^{***}( x^{**}
, y_{\beta}), z^{*}\rangle\\&=&\lim_{\beta}\langle x^{**}, f^{**}(y_{\beta},  z^*)\rangle\\&=&\lim_{\beta}\lim_{\alpha}\langle f^{**}(y_{\beta}, z^*),  x_{\alpha}\rangle\\&=&\lim_{\beta}\lim_{\alpha}
\langle z^{*}, f( x_{\alpha}, y_{\beta})\rangle\\&=&\langle f^{r***r}( x^{**}, y^{**}) , z^{*}\rangle.\eea
It follows that $f$ is regular and this completes the proof.
\end{proof}
As the first application of Theorem ~\ref{T1}, we may present the following results of \cite {DPV}, with a direct proof.

\begin{corollary}  [{\cite[Propositions 4.1 and  4.4]{DPV}}]\label{C1}
For a bounded bilinear map $f :X \times Y\longrightarrow Z$, the following statements are equivalent:

$(i)$  $f$ and $f^{r*}$ are regular.

$(ii)$ $f^{***r*r}=f^{r*r***}$.

$(iii)$ $ f^{****} (Z^{***},  X^{**}) \subseteq Y^{*}$.
\end{corollary}
\begin{proof} The implication  $(i)\Rightarrow (ii)$ follows  trivially.

If $(ii)$ holds then  $f^{****}=f^{r*r**}$. Indeed, for every  $x^{**}\in X^{**}, y^{**}\in Y^{**}$ and $z^{***}\in Z^{***},$
\bea \langle f^{****}(z^{***}, x^{**}), y^{**}\rangle&=&\langle f^{***r*r}(y^{**}, z^{***}), x^{**}\rangle\\&=&\langle f^{r*r***}(y^{**}, z^{***}), x^{**}\rangle\\&=&\langle f^{r*r**}(z^{***}, x^{**}), y^{**}\rangle.\eea
As $f^{r*r**}(Z^{***}, X^{**})$ always lies in $Y^*$, we have reached  $(iii)$.

 For $(iii)\Rightarrow (i)$, since $f^{****} (Z^{*}, X^{**}) \subseteq f^{****} (Z^{***},  X^{**}) \subseteq Y^{*}$, Theorem~\ref{T1} implies the regularity of $f$, or equivalently $f^{r*****}=f^{****r},$ from which

 \[(f^{r*})^{****}(X^{**},  Z^{***})= f^{r*****}(X^{**}, Z^{***})=f^{****} (Z^{***}, X^{**})\subseteq Y^*.\] Now the regularity of $f^{r*}$ follows trivially again by   Theorem~\ref{T1}.
\end{proof}

 The main result of   Arikan's paper, (\cite[Theorem 2]{AR2}),  is a criterion for the Arens regularity of bilinear mappings which applied  it to give  a train of results on the regularity of the algebra $l^1$ with pointwise multiplication, the algebra $L^{\infty}(G)$ with convolution, where $G$ is a compact group, and the trace-class algebra. Now, we give  her result as a consequence of  Corollary~\ref{C1}. It is worthwhile discussing that the assumption  that she applied in \cite[Theorem 2]{AR2} has actually more byproducts; indeed, as we shall see  in the next corollary,  both $f$ and $f^{r*}$ are regular, however she deduces  merely the regularity of $f$.

\begin{corollary}  [see {\cite[Theorem 2]{AR2}}]\label{C2}
Let $f :X \times Y\longrightarrow Z$ and $g :X \times W\longrightarrow Z$ be bounded bilinear mappings and let $h :Y\longrightarrow W$ be a bounded linear mapping such that $f(x, y)=g(x, h(y))$, for all $x\in X, y\in Y$. If $h$ is weakly compact, then both $f$ and $f^{r*}$ are regular.
\end{corollary}
\begin{proof}
Using the equality $f(x, y)=g(x, h(y))$, a standard argument applies to show that $f^{****}=h^{***}\circ g^{****}$. The weak compactness of $h$ implies that of $h^*$, from which we  have $h^{***}(W^{***})\subseteq Y^*$. Therefore
\[f^{****} (Z^{***},  X^{**})=h^{***}(g^{****} (Z^{***},  X^{**}))=h^{***}(W^{***})\subseteq Y^{*}.\]
Now Corollary~\ref{C1} implies that both $f$ and $f^{r*}$ are regular, as claimed.
\end{proof}

\section{Arens regularity of  module actions}
Let $A$ be a Banach algebra, $X$ be a Banach space and let  $\pi_{1} : A\times
X\longrightarrow X$  be a bounded bilinear map. Then the pair $(\pi_1, X)$  is said to be a left  Banach $A-$module if  $\pi_1(ab, x)=\pi_1(a, \pi_1(b, x)),$ for every $a, b\in A, x\in X.$
A right Banach $A-$module $(X, \pi_2)$ may be defined similarly.
A triple $(\pi_1, X, \pi_2)$ is said to be a Banach $A-$module if $(\pi_1, X)$ and  $(X, \pi_2)$
are left and  right Banach $A-$modules, respectively, and for every $a, b\in A, x\in X,$
\[\pi_1(a, \pi_2(x, b))=\pi_2(\pi_1(a, x), b).\]

If so, then trivially $({\pi_2}^{r*r}, X^*, \pi_1^*)$ is the dual Banach $A-$module of $(\pi_1, X, \pi_2)$. Also
 $(\pi_{1}^{***}, X^{**}, \pi_{2}^{***})$  is a Banach $(A^{**},\Box)-$module as well as $(\pi_{1}^{r***r}, X^{**}, \pi_{2}^{r***r})$  is a Banach  $(A^{**},\lozenge)-$module. To verify these, one may carefully check that the various associativity rules follow from the analogous rules that show that $(\pi_1, X, \pi_2)$ is a Banach $A-$module.

Now if we continue  dualizing we shall reach  $(\pi_{2}^{***r*r}, X^{***}, \pi_{1}^{****})$  as the dual Banach $(A^{**}, \Box)-$module of  $(\pi_{1}^{***}, X^{**}, \pi_{2}^{***})$ as well as $(\pi_{2}^{r****r}, X^{***}, \pi_{1}^{r***r*})$ as the dual Banach $(A^{**},\lozenge)-$module of $(\pi_{1}^{r***r}, X^{**}, \pi_{2}^{r***r})$.

  From now on, for a Banach $A-$module $X$ we usually use the above mentioned  canonical module actions on $X^*$, $X^{**}$ and $X^{***}$, unless otherwise stipulated  explicitly. However, for  brevity of notation,  when there is no risk of confusion, we may use them without the specified module actions.\\

 According to the inclusions
\[\pi_{2}^{***r*r}(A^{**}, X^{*})=\pi_{2}^{**}(A^{**}, X^{*})\subseteq X^*\] and
\[\pi_{1}^{r***r*}(X^{*}, A^{**})=\pi_{1}^{r**r}(X^*, A^{**})\subseteq X^*,\]
we deduce that $(\pi_2^{**}, X^*)$ is a left  $(A^{**}, \Box)-$submodule of $X^{***}$,
while  $(X^*, \pi_1^{r**r})$ is a  right $(A^{**}, \lozenge)-$submodule of $X^{***}$.  The next result determines when   $(\pi_2^{**}, X^*, \pi_1^{r**r})$  is actually a  $(A^{**}, \Box)$ (respectively,  $(A^{**}, \lozenge)$)$-$submodule of $X^{***}$.
\begin{proposition}\label{p1} Let $(\pi_1, X, \pi_2)$ be a Banach $A-$module. Then

 $(i)$ $X^{*}$ is a   $(A^{**}, \Box)-$submodule of $X^{***}$ if and only if  $\pi_{1}$ is regular.

 (ii) $X^{*}$ is a  $(A^{**}, \lozenge)-$submodule of $X^{***}$ if and only if $\pi_{2}$ is regular.
\end{proposition}
\begin{proof}
 We prove only $(i)$, the other one uses the same argument. It is trivial that $X^{*}$ is a right  $(A^{**}, \Box)-$submodule of $X^{***}$ if and only if $ \pi_{1}^{****}(X^*, A^{**}) \subseteq X^*$ and by  Theorem~\ref{T1} this is nothing but the regularity of $\pi_{1}$; which establishes $(i)$.
\end{proof}

 As a rapid consequence  of the latter result, we examine  it for  $\pi_1=\pi$ and $\pi_2=\pi^r$ on $X=A$ and we thus have the following result of Dales, Rodriguez-Palacios and Velasco, \cite
 {DPV}.

\begin{corollary} [{\cite[Proposition 5.2]{DPV}}] For a  Banach algebra $A$ the following statements are equivalent:

(i) $A$ is  Arens regular.

(ii) $A^*$ is a  $(A^{**}, \Box)-$submodule of $A^{***}$.

(iii) $A^*$ is a $(A^{**}, \lozenge)-$submodule of $A^{***}$.

\end{corollary}

Let $(\pi_1, X, \pi_2)$ be a Banach $A-$module. As we  mentioned just before  Proposition~\ref{p1}, $(\pi_2^{**}, X^*)$ and $(X^*, \pi_1^{r**r})$ are  left  $(A^{**}, \Box)-$submodule and  right $(A^{**}, \lozenge)-$submodule of $X^{***}$, respectively. So if we assume that $A$ is Arens regular then $(\pi_2^{**}, X^*)$ and  $(X^*, \pi_1^{r**r})$ are left and right $A^{**}-$modules, respectively. The next result deals with  the question when   $(\pi_2^{**}, X^*, \pi_1^{r**r})$ is actually a $A^{**}-$module.
\begin{proposition}\label{p2}
Let $A$ be  Arens regular and let $(\pi_1, X, \pi_2)$ be a Banach $A-$module. Then $(\pi_2^{**}, X^*, \pi_1^{r**r})$ is a Banach $A^{**}-$module if and only if the  bilinear map \[\theta_x :A\times A\longrightarrow  X,\ \ \ \theta_x(a, b)=\pi_1( a, \pi_2(x, b))=\pi_2(\pi_1(a, x), b) \ \ \ \ (a, b\in A) \] is  regular,  for all $x\in X$.
\end{proposition}
\begin{proof}
 The triple $(\pi_2^{**}, X^*, \pi_1^{r**r})$ is a Banach $A^{**}-$module  if and only if  for all $a^{**}, b^{**}\in A^{**}, x^*\in X^*$ and $x\in X$,  \[\langle{\pi_1}^{r**r}({\pi_2}^{**}(b^{**}, x^*), a^{**}), x\rangle=\langle{\pi_2}^{**}(b^{**}, {\pi_1}^{r**r}(x^*, a^{**})), x\rangle.\]
Let $\{a_\alpha\}$ and $\{b_\beta\}$ be two nets in $A$ that converge to $a^{**}$ and $b^{**}$, respectively, in  the  $w^*-$topology of $A^{**}$. Then a direct verification reveals that\\
 \[\langle{\pi_1}^{r**r}({\pi_2}^{**}(b^{**}, x^*), a^{**}), x\rangle=\lim_\alpha\lim_\beta\langle x^*, \pi_2(\pi_1(a_\alpha, x), b_\beta)\rangle=\lim_\alpha\lim_\beta\langle x^*, \theta_x(a_\alpha, b_\beta)\rangle\] and
\[\langle{\pi_2}^{**}(b^{**}, {\pi_1}^{r**r}(x^*, a^{**})), x\rangle=\lim_\beta\lim_\alpha\langle x^*, \pi_1( a_\alpha, \pi_2(x, b_\beta))=\lim_\beta\lim_\alpha\langle x^*, \theta_x(a_\alpha, b_\beta)\rangle.\]
Thus   $(\pi_2^{**}, X^*, \pi_1^{r**r})$ is a Banach $A^{**}-$module if and only if $\theta_x$ is regular, for all $x\in X$.
\end{proof}

As a consequence, we give the next result of Bunce and Paschke, \cite {BP}.
\begin{corollary} [{\cite[Proposition 1.1]{BP}}]
Let $A$ be a $C^*-$algebra  and let $(\pi_1, X, \pi_2)$ be a Banach $A-$module. Then $(\pi_2^{**}, X^*, \pi_1^{r**r})$ is a Banach $A^{**}-$module.
\end{corollary}
\begin{proof}
By Proposition~\ref{p2} it is enough to show that $\theta_x :A\times A\longrightarrow X$ is regular, for all $x\in X$. Applying  Theorem~\ref{T1}, the regularity of $\theta_x$ is equivalent to the weak compactness of the linear mapping $a\mapsto {\theta_x}^*(x^*, a) :A\longrightarrow A^*$, which is guaranteed by the fact that, every bounded
linear map from a $C^{*}-$algebra to the predual of a $W^{*}-$
algebra is automatically weakly compact (see \cite[Corollary II 9]{AK}).
\end{proof}
\begin{remarks}$\ \ \ $

$(i)$ If either $\pi_1$ or  $\pi_2$ is regular then trivially $\theta_x$ is regular, for each $x\in X$. Therefore, in the Arens regular setting for $A$, using  Proposition~\ref{p2},  $X^*$ is a $A^{**}-$module, which is actually a $A^{**}-$submodule of $X^{***}$, (see Proposition~\ref{p1}).

$(ii)$  In the Arens regular setting for $A$, if  $\theta_x$ is regular for each $x\in X$  then  $(\pi_2^{**}, X^*, \pi_1^{r**r})$ is a Banach $A^{**}-$module by  Proposition~\ref{p2}, which gives rise naturally the  dual Banach $A^{**}-$module  $(\pi_1^{r***r}, X^{**}, \pi_2^{***})$. On the other hand we have the canonical Banach $A^{**}-$modules  $(\pi_1^{***}, X^{**}, \pi_2^{***})$ and $(\pi_1^{r***r}, X^{**}, \pi_2^{r***r})$. It is worthwhile to note that the involved dual module  does not coincide with  the canonical  modules, in general. However, it is obvious that it  coincides with one of them  if (and only if) either $\pi_1$ or $\pi_2$ is  regular.
It is shown in \cite[Proposition 1.2]{BP} that the involved module actions  coincide when $A$ is a $C^*-$algebra and $X^*$ is weakly sequentially complete; however, according to what we shall demonstrate in the following, under these conditions  on $A$ and $X$ the module actions $\pi_1$ and $\pi_2$ are regular.

$(iii)$ Let $(\pi_1, X, \pi_2)$ be a Banach $A-$module, if $A$ is a $C^*-$algebra and $X^*$ is weakly sequentially complete, then by \cite[Theorem 4.2]{ADG} the bounded linear mappings $a\longmapsto \pi_1^*(x^*, a)$ and  $a\longmapsto \pi_2^{r*}(x^*, a)$ from $A$ to $X^*$ are  weakly compact for all $x^{*} \in X^{*}$. Now  applying  Theorem~\ref{T1} for $\pi_1$ and $\pi_2^r$ shows that $\pi_1$ and $\pi_2$ are regular. Note that the weak sequential completeness of $X^*$ is essential and can not be removed in general. For instance, let $A$ be the $C^*-$algebra of compact operators on a separable, infinite dimensional Hilbert space $H$ and let $X$ be the trace-class operators on $H$; then a direct verification reveals that the usual  $A-$module action  on $X$ is not regular (see the example just before  Proposition 2.1 in \cite{BP}).
\end{remarks}

Dales, Rodriguez-Palacios and Velasco in \cite[Proposition 4.5]{DPV} (see also, \cite[Theorem 4]{AR1} and \cite[Theorem 3.1]{U}) have shown that, if $A$ is Arens regular with a bounded left approximate identity, then $\pi^{r*} :A^*\times A\longrightarrow A^*$ is regular if and only if $A$ is reflexive as a Banach space. (Note that,  the equality $ \pi^{***r*r}=\pi^{r*r***}$ which they used in \cite[Proposition 4.5]{DPV} is equivalent to the regularity of  both $A$ and $\pi^{r*}$, see Corollary~\ref{C1}). In the next proposition we generalize their result, which also shows that the hypothesis of Arens regularity of $A$ in \cite[Proposition 4.5]{DPV} is superfluous.  Before proceeding, we  recall that, when $(\pi_1, X)$ is a left Banach $A-$module then a bounded net $\{e_\alpha\}$ in $A$ is said to be a left approximate identity for $X$, if $\pi_{1}(e_\alpha, x)\longrightarrow x$, for each $x\in X$. As a consequence of the so-called Cohen Factorization Theorem, see \cite D, it is  known that a bounded left approximate identity of $A$ is that of $X$ if and only if  $\pi_{1}(A, X)=X$. The same situation  happens for the right Banach $A-$module $(X, \pi_2)$.

\begin{proposition} \label{p3} Let $(\pi_{1}, X)$ and $(X, \pi_{2})$ be  left and  right   Banach  $A-$modules respectively.

 $(i)$ If $A$ has a bounded left approximate identity for $X$, then $\pi_{1}^{r*}$ is regular
if and only if  $X$ is reflexive.

 $(ii)$ If $A$ has a bounded right approximate identity for $X$, then $\pi_{2}^{*}$ is regular
if and only if  $X$ is reflexive.

\end{proposition}

\begin{proof} $(i)$ If $X$ is reflexive, then \[(\pi_{1}^{r*})^{****}(A^{**}, X^{***})=\pi_{1}^{r*****}(A^{**}, X^{*})=\pi_{1}^{r**}(A^{**}, X^{*})\subseteq X^*;\] now Theorem~\ref{T1} implies that $\pi_{1}^{r*}$ is regular. For the converse, first note that,
 if $\{e_{\alpha}\}$ is a  bounded
 left approximate identity for $X$ (in $A$), then a direct verification reveals that for every $x^{**}\in X^{**}$, $\pi_{1}^{r***}(x^{**}, e^{**})=x^{**}$,  in which $e^{**}$ is a $w^*-$cluster point of $\{e_{\alpha}\}$ in $A^{**}$. The latter identity in turn implies that $x^{***}=\pi_1^{r*****}(e^{**}, x^{***})$, for every $x^{***}\in X^{***}$. Now   the regularity of $\pi_{1}^{r*}$ (again by Theorem~\ref{T1}) shows that  $\pi_1^{r*****}(e^{**}, x^{***})\in X^*$. We thus have $x^{***}\in X^*$,  or equivalently  $X$ is reflexive. The proof of $(ii)$ is very similar to that of $(i)$.
\end{proof}

\section{The second adjoint of a derivation}
Let $(\pi_1, X, \pi_2)$ be a Banach $A-$module. A bounded linear mapping $D :
A\longrightarrow X^*$ is said to be a derivation if  $ D(ab)=D(a)\cdot b + a\cdot D(b)$, for each $a,b\in A,$
 or equivalently \[D(ab) = \pi_1^*(D(a),b) + \pi_2^{r*r}(a,D(b)).\]

In this section we deal with the question of when the second adjoint $D^{**} :A^{**}\longrightarrow X^{***}$ of  $D : A\longrightarrow X^{*}$ is again a derivation. We recall from the beginning of the last section that  $(\pi_{2}^{***r*r}, X^{***}, \pi_{1}^{****})$ and  $(\pi_{2}^{r****r}, X^{***}, \pi_{1}^{r***r*})$ are  our canonical Banach $(A^{**}, \Box)-$module and  $(A^{**},\lozenge)-$ module, respectively. Hence, $D^{**} :(A^{**}, \Box)\longrightarrow X^{***}$ is a derivation if and only if for every $a^{**}, b^{**}\in A^{**},$
\[ D^{**}(a^{**}\Box \ b^{**})=
\pi_{1}^{****}(D^{**}(a^{**}),b^{**}) +
\pi_{2}^{***r*r}(a^{**},D^{**}(b^{**})).\]
Similarly,  $D^{**} :(A^{**}, \lozenge)\longrightarrow X^{***}$ is a derivation if and only if
\[ D^{**}(a^{**}\lozenge \ b^{**})=
\pi_{1}^{r***r*}(D^{**}(a^{**}),b^{**}) +
\pi_{2}^{r****r}(a^{**},D^{**}(b^{**})).\]

We commence with the following lemma.
\begin{lemma}\label{l} Let $(\pi_1, X, \pi_2)$ be a Banach $A-$module, and let $D :
A\longrightarrow X^{*}$ be a derivation. Then for all $a^{**}, b^{**}\in A^{**},$

$(i)$ $D^{**}(a^{**}\Box \ b^{**})=\pi_{1}^{****}( D^{**}(a^{**}), b^{**})+\pi_{2}^{r*r***}(a^{**},  D^{**}(b^{**}))$ and

$(ii)$ $D^{**}(a^{**}\lozenge \ b^{**})=\pi_{1}^{*r***r}( D^{**}(a^{**}), b^{**})+\pi_{2}^{r****r}(a^{**},  D^{**}(b^{**}))$.

\end{lemma}
\begin{proof}
$(i)$ Let $\{a_\alpha\}$ and $\{b_\beta\}$ be two nets in $A$ that converge  to $a^{**}$ and $b^{**}$ in the  $w^*-$topology of $A^{**}$, respectively. Then for each $x^{**}\in X^{**}$ we have
\bea \langle D^{**}(a^{**}\Box \ b^{**}), x^{**}\rangle &=&\lim_\alpha\lim_\beta\langle D(a_\alpha b_\beta), x^{**}\rangle\\
                                         &=&\lim_\alpha\lim_\beta\langle\pi_{1}^{*}(D(a_\alpha),b_\beta)+\pi_{2}^{r*r}(a_\alpha,D(b_\beta)), x^{**}\rangle\\
                                         &=&\lim_\alpha\lim_\beta\langle\pi_{1}^{**}(x^{**}, D(a_\alpha)) , b_\beta\rangle+\lim_\alpha\lim_\beta\langle\pi_{2}^{r*r*}(x^{**}, a_\alpha), D(b_\beta)\rangle\\
                                         &=&\lim_\alpha\langle\pi_{1}^{**}(x^{**}, D(a_\alpha)) , b^{**}\rangle+\lim_\alpha\langle\pi_{2}^{r*r*}(x^{**}, a_\alpha), D^{**}(b^{**})\rangle\\
                                         &=&\langle\pi_{1}^{***}(b^{**}, x^{**}) , D^{**}(a^{**})\rangle+\langle\pi_{2}^{r*r**}( D^{**}(b^{**}), x^{**}), a^{**}\rangle\\
                                         &=&\langle\pi_{1}^{****}( D^{**}(a^{**}), b^{**})+\pi_{2}^{r*r***}(a^{**},  D^{**}(b^{**})), \ x^{**}\rangle.\eea
  Hence $D^{**}(a^{**}\Box \ b^{**})=\pi_{1}^{****}( D^{**}(a^{**}), b^{**})+\pi_{2}^{r*r***}(a^{**}, D^{**}(b^{**}))$. A similar argument  applies for $(ii)$.
\end{proof}

Dales, Rodriguez-Palacios and Velasco, in the main theorem  of their paper,  \cite[Theorem 7.1]{DPV}, have shown that the second adjoint $D^{**} :(A^{**}, \Box) \longrightarrow A^{***}$ of  the derivation $D : A\longrightarrow A^{*}$ is  a derivation
if and only if $D^{**}(A^{**})\cdot A^{**}\subseteq A^*$, or equivalently $\pi^{r****}(D^{**}(A^{**}),  A^{**})\subseteq A^*$, where  $\pi$ is the multiplication of $A$. The next theorem extends it with a direct proof.
\begin{theorem}\label{T2}
Let $(\pi_1, X, \pi_2)$ be a Banach $A-$module and let $D :
A\longrightarrow X^{*}$ be a derivation.

$(i)$ $D^{**} :(A^{**}, \Box)\longrightarrow X^{***}$ is a derivation if and only if $ \pi_{2}^{****}(D^{**}(A^{**}),  X^{**})\subseteq A^*$.

$(ii)$ $D^{**} :(A^{**}, \lozenge)\longrightarrow X^{***}$ is a derivation if and only if $ \pi_{1}^{r****}(D^{**}(A^{**}), X^{**})\subseteq A^*$.
\end{theorem}
\begin{proof}
 We  only prove $(i)$. Let  $a^{**}, b^{**}\in A^{**}$ and $x^{**}\in X^{**}$ be arbitrary. Applying Lemma~\ref{l} and the equations just before  it, one may  deduce that $D^{**} :(A^{**}, \Box)\longrightarrow X^{***}$  is a derivation if and only if
\[\langle\pi_{2}^{***r*r}(a^{**},  D^{**}(b^{**})), x^{**}\rangle=\langle\pi_{2}^{r*r***}(a^{**},D^{**}(b^{**})), x^{**}\rangle;\]
which  holds  if and only if
\[\langle\pi_{2}^{****}( D^{**}(b^{**}), x^{**}), a^{**}\rangle=\langle\pi_{2}^{r*r**}( D^{**}(b^{**}), x^{**}), a^{**}\rangle,\] or equivalently, \[\pi_{2}^{****}( D^{**}(b^{**}), x^{**})=\pi_{2}^{r*r**}( D^{**}(b^{**}), x^{**}).\] It should be noted that $\pi_{2}^{****}( D^{**}(b^{**}), x^{**})\in A^{***}$, while $\pi_{2}^{r*r**}( D^{**}(b^{**}), x^{**})\in A^*$ and also
$\pi_{2}^{****}( D^{**}(b^{**}), x^{**})_{|_A}=\pi_{2}^{r*r**}( D^{**}(b^{**}), x^{**}).$ Thus $D^{**}$ is a derivation if and only if $\pi_{2}^{****}( D^{**}(b^{**}), x^{**})\in A^{*}$.\\
\end{proof}

As  immediate consequences of the  Theorem~\ref{T2} we have the next corollaries.

\begin{corollary}Let $(\pi_1, X, \pi_2)$ be a Banach $A-$module and let $D :
A\longrightarrow X^{*}$ be a derivation. If $A$ is Arens regular then the following statements are equivalent:

$(i)$  $D^{**} :A^{**}\longrightarrow X^{***}$ is a derivation.

$(ii)$ $ \pi_{2}^{****}(D^{**}(A^{**}),  X^{**})\subseteq A^*$.

$(iii)$ $\pi_{1}^{r****}(D^{**}(A^{**}),  X^{**})\subseteq A^*$.
\end{corollary}
\begin{corollary} Let $(\pi_1, X, \pi_2)$ be a Banach $A-$module, and let $D :
A\longrightarrow X^{*}$ be a derivation.

$(i)$ If both $\pi_{2}$ and $\pi_{2}^{r*}$ are Arens regular then $D^{**} :(A^{**}, \Box)\longrightarrow X^{***}$ is a derivation.

$(ii)$ If both $\pi_{1}$ and $\pi_{1}^{*}$ are Arens regular then $D^{**} :(A^{**}, \lozenge)\longrightarrow X^{***}$ is a derivation.
\end{corollary}
\begin{proof}
$(i)$ If both $\pi_{2}$ and $\pi_{2}^{r*}$ are regular then $\pi_{2}^{****}(X^{***}, X^{**})\subseteq A^*$ by Corollary~\ref{C1}. In particular, $\pi_{2}^{****}(D^{**}(A^{**}), X^{**})\subseteq A^*$ and the conclusion  follows from Theorem~\ref{T2}. A similar proof applies for $(ii)$.

\end{proof}
However the  regularity of both $\pi_{2}$ and $\pi_{2}^{r*}$ in $(i)$ (in turn $\pi_{1}$ and $\pi_{1}^{*}$ in $(ii)$) of the latter corollary ensure that $D^{**}$ is a derivation, but it seems that these conditions impose  rather a strong requirement on the module actions. As another application of  Theorem~\ref{T2}, we give the following generalization of  \cite[Corollary 7.2(i)]{DPV} which share the imposed requirements on both module action and the derivation, as well.

\begin{corollary} Let $(\pi_1, X, \pi_2)$ be a Banach $A-$module and let $D :
A\longrightarrow X^{*}$ be a weakly compact derivation.

$(i)$ If $\pi_{2}$ is regular then $D^{**} :(A^{**}, \Box)\longrightarrow X^{***}$ is a derivation.

$(ii)$ If $\pi_{1}$ is  regular then $D^{**} :(A^{**}, \lozenge)\longrightarrow X^{***}$ is a derivation.
\end{corollary}
\begin{proof} $(i)$ Since $D :A\longrightarrow X^{*}$ is weakly compact, $D^{**}(A^{**})\subseteq X^*$. On the other hand  the regularity of  $\pi_{2}$ implies that
$\pi_{2}^{****}(X^{*}, X^{**})\subseteq A^*$ (see Theorem~\ref{T1}). Therefore  $\pi_{2}^{****}(D^{**}(A^{**}),  X^{**})\subseteq A^*$, which means that $D^{**} :(A^{**}, \Box)\longrightarrow X^{***}$ is a derivation.
\end{proof}
 We conclude this section with some observations on inner derivations. A linear mapping  $D :
A\longrightarrow X^{*}$ is said to be an inner derivation if  $D(a) = x^{*}\cdot a - a\cdot x^{*}$ for some
$x^{*}\in X^{*}$,  or equivalently $D(a) = \pi_{1}^{*}(x^{*},a) -
\pi_{2}^{r*r}(a,x^{*})$, for all $a \in A$. Then it follows trivially that for every $a^{**}\in A^{**},$
\[D^{**}(a^{**})=\pi_{1}^{****}(x^{*},a^{**})-\pi_{2}^{r****}(x^{*}, a^{**}).\]
If, in addition, we assume that $\pi_{1}$ and $\pi_{2}$ are regular then Theorem~\ref{T1} implies that the right hand side of the latter equality  belongs to $X^*$;  thus  $D^{**}(A^{**})\subseteq X^{*}$, or equivalently, $D$ is weakly compact. Furthermore, in the same situation the latter identity together with the equalities $\pi_{1}^{****}=\pi_{1}^{r***r*}$ and $\pi_{2}^{r****r}=\pi_{2}^{***r*r}$ imply that,
\[D^{**}(a^{**})=
\pi_{1}^{****}(x^{*},a^{**})-\pi_{2}^{***r*r}(a^{**},x^{*})\]
and also
\[D^{**}(a^{**})=
\pi_{1}^{r***r*}(x^{*},a^{**})-\pi_{2}^{r****r}(a^{**},x^{*});\]
which means that both $D^{**} :(A^{**}, \Box)\longrightarrow X^{***}$ and  $D^{**} :(A^{**}, \lozenge)\longrightarrow X^{***}$ are also  inner derivations. We summarize these observations in the next result which is a generalization of \cite[Proposition 6.1]{DPV}.
\begin{proposition} Let $X$ be a $A-$module whose module actions are regular. Then every inner derivation
$D :A\longrightarrow X^{*}$ is weakly compact; moreover,  $D^{**} :(A^{**}, \Box)\longrightarrow X^{***}$ and  $D^{**} :(A^{**}, \lozenge)\longrightarrow X^{***}$ are also inner derivations.
\end{proposition}

\section*{Acknowledgments}This paper was completed while the second author was visiting the University of New South Wales in Australia. He would like to thank the School of Mathematics and Statistics, especially Professor Michael Cowling for their generous hospitality. The valuable discussion of Professor Garth Dales in his short visit to Mashhad University is also acknowledged.

\end{document}